# ON $H^1$ OF FINITE DIMENSIONAL ALGEBRAS

Claude Cibils

### Abstract

We review results on the first Hochschild cohomology vector space of a finite dimensional algebra, in particular for path algebras modulo a "pre-generated" ideal. In case of a monomial algebra whose quiver has no oriented cycles, a dimension formula is provided. The general monomial case is considered in a paper with M. Saorin [9]

## 1 Introduction

Let $k$ be a field and let $\Lambda$ be a finite dimensional $k$-algebra. Recently several authors have provided interest on the first Hochschild cohomology group of $\Lambda$, see for instance [1, 2, 13, 20, 21]. One aim of this note is to present an account of known results on $H^1$ which have been somehow screened by more general statements. We also present a new dimension formula for $H^1$ of monomial algebras in case the quiver has no oriented cycles, and intend to provide an introduction to the subject.

Note that important and recent results concerning the $H^1$ of a finite dimensional algebra can be founded in [20, 3]. In particular for algebras whose quiver have no oriented cycles, the dimension of $H^1$ is bounded below by the rank of the free component of the abelianisation of the fundamental group, see [20]. We also mention that very interesting examples have been produced by Buchweitz and Liu [3]: they show that there exist algebras having zero $H^1$ and loops in their quiver.

First we recall that a $\Lambda$-bimodule $X$ is a left and right $\Lambda$-module verifying $(\lambda x)\mu = \lambda(x\mu)$ for every $\lambda, \mu$ in $\Lambda$ and $x$ in $X$, such that the action of $k$ is central: $ax = xa$ for $a \in k$ and $x \in X$. Note that if $\Lambda^e = \Lambda \otimes_k \Lambda^{op}$ is the enveloping algebra, $\Lambda^e$-left modules and $\Lambda$-bimodules are the same.

The Hochschild cohomology vector spaces are extensions vector spaces (see for instance [10] or [19]):

$$H^i(\Lambda, X) = \mathsf{Ext}^i_{\Lambda^e}(\Lambda, X).$$

There is a standard projective resolution of $\Lambda$ as a $\Lambda^e$-module which provides the usual Hochschild complex of cochains computing the Hochschild cohomology $H^i(\Lambda, X)$, see [10, 18, 19].

At zero-degree we have

$$H^0(\Lambda, X) = \mathsf{Hom}_{\Lambda^e}(\Lambda, X) = X^\Lambda$$

where $X^\Lambda = \{x \in X | \lambda x = x\lambda \ \forall \lambda \in \Lambda\}$. Indeed, to $\varphi \in \mathsf{Hom}_{\Lambda^e}(\Lambda, X)$ we associate $\varphi(1)$ which belongs to $X^\Lambda$. In particular if the $\Lambda$-bimodule $X$ is the algebra itself, $H^0(\Lambda, \Lambda)$ is the center $\mathbf{Z}(\Lambda)$ of $\Lambda$.

The trivial case corresponds to an algebra $\Lambda$ which is projective as a $\Lambda^e$-left module. This occurs only when $\Lambda$ is semisimple (a product of matrix algebras over finite dimensional skew-fields) and the centers of the skew-fields are separable extensions of $k$. Such algebras are called separable. A zero-length projective resolution of a separable algebra $\Lambda$ shows that $H^i(\Lambda, X) = 0$ for $i \geq 1$. See [10] for proofs.

The first Hochschild cohomology vector space $H^1(\Lambda, X)$ is the quotient of the derivations by the inner ones: a derivation $f : \Lambda \to X$ is a linear map verifying $f(\lambda\mu) = \lambda f(\mu) + f(\lambda)\mu$ for $\lambda$ and $\mu$ in $\Lambda$. It is inner if there exists a $x \in X$ such that $f(\lambda) = \lambda x - x\lambda$.

We also record that $H^2(\Lambda, \Lambda)$ is related with the deformation theory of $\Lambda$, see [14].

In the following sections we will study $H^1$ for algebras of the form $kQ/I$ where $Q$ is a quiver, $kQ$ is its path algebra and $I$ is an ideal of $kQ$. Actually by an observation of P. Gabriel [12] any finite dimensional $k$-algebra over an algebraically closed field is Morita equivalent to an algebra of this sort. The interesting point is that the ideal $I$ can be chosen to be admissible, that is contained in $F^2$ where $F$ is the two-sided ideal generated by the arrows of $Q$. We recall these well known results in the next section.

We consider the case where $I$ is a "pre-generated" ideal, the definition is given at section 3. This includes the cases $I = 0$ whenever $Q$ has no oriented cycles, any ideal of a narrow quiver, and some other cases. An explicit dimension formula for $H^1(kQ/I, kQ/I)$ can be given, note that these results belongs to [4]. D. Happel considered in [17] monomial Schurian almost commutative algebras which corresponds actually to instances of pre-generated ideals.

In the fifth section we study monomial algebras, that is path algebras of a quiver $Q$ modulo a two-sided ideal $I$ admitting a generating set of paths of length at least two. When $Q$ has no oriented cycles, we show that a formula computing $\dim_k H^1(kQ/I, kQ/I)$ follows immediately from the considerations made in [5] where rigid monomial algebras are classified through a careful study of $H^2$. In a paper with Manuel Saorin [9] we provide the general formula computing the dimension of the first Hochschild cohomology vector space of a monomial algebra whith an arbitrary quiver. In the sixth section we specialize to truncated algebras the results on monomial algebras.

In the last section we consider finite posets and their incidence algebras: it is known that the Hochschild cohomology is isomorphic to the usual cohomology of the simplicial complex provided by the chains (see [6, 15, 16]).

In section 4 we consider a tensor algebra $T$ over a separable algebra $E$. Actually path algebras of any finite quiver are of this sort. We study the case of a $T$-bimodule of coefficients $X$ which is finitely generated as an $E$-bimodule. A formula for the dimension of $H^1(T, X)$ can be obtained directly from [5]. The important point is that $T$ has a length one projective resolution as a $T$-bimodule. These results are useful when considering monomial algebras.

Acknowledgements: I thank Manuel Saorin for expressing his interest on $H^1$ of monomial algebras. I also thank Belkacem Bendiffalah for discussions concerning the cohomology of categories.

## 2   Quivers, ideals and categories

This section is an account of well known results which are usually implicit in most research papers, see also [20].

A quiver is given by two sets $Q_0$, $Q_1$ and two maps $s, t : Q_1 \to Q_0$. An arrow $a \in Q_1$ of the corresponding oriented graph has source vertex $s(a)$ and terminal vertex $t(a)$. A path $\gamma$ of $Q$ is a finite sequence of concatenated arrows. Each path has a source vertex $s(\gamma)$ and a terminal vertex $t(\gamma)$. An oriented cycle is a path which have same source and terminal, and we agree that a vertex is a length zero oriented cycle.

The quiver algebra is the vector space of finite linear combinations of paths. The multiplication of two paths is their concatenation if it can be done and 0 otherwise. Each vertex is an idempotent, the sum of all the vertices is the unit element in case $Q_0$ is finite. Some authors are interested in the categorical point of view, we provide an overview of this aspect and its relation with finite dimensional algebras.

Let $C$ be a small $k$-category (morphisms are $k$-vector spaces and the composition is $k$-bilinear). We denote Ob$C$ the set of objects and $_yC_x$ the vector space of morphisms from $x$ to $y$. The Gabriel algebra [11] of $C$ is $\bigoplus_{x,y \in \text{Ob}C} {_yC_x}$. The product corresponds to the matrix product: let $f = (_yf_x)$ and $g = (_yg_x)$ be elements of Gabriel's algebra. The $(y, x)$ component of their product is $\sum_{v \in \text{Ob}C} {_yg_u}\ {_uf_x}$. If Ob$C$ is a finite set the Gabriel algebra has a unit element given by the identities.

A quiver $Q$ provides a free $k$-category $C_Q$ as follows: the objects of $C_Q$ are the vertices of $Q$, and $_y(C_Q)_x$ is the vector space of finite linear combinations of paths from $x$ to $y$. Concatenation of paths provides the composition of $C_Q$ on the basis and is extended $k$-linearly. Note that $C_Q$ is free in the sense that any $k$-linear functor $F : C_Q \to D$ is completely determined by its values on the objects and on the arrows, provided $F(a) \in {_{F(t(a))}D_{F(s(a))}}$ for each $a \in Q_1$.

Of course, the Gabriel algebra of $C_Q$ is exactly the path algebra $kQ$ that we have defined directly before.

In [12] P. Gabriel noted that if $k$ is algebraically closed and if the Gabriel algebra of a $k$-category $C$ is finite dimensional, then a skeleton of $C$ is a quotient of a free category $C_Q$ (a skeleton is given by a choice of one object in each isomorphism class of objects and considering the corresponding full subcategory).

In other words, let $\Lambda$ be a finite dimensional algebra over an arbitrary field $k$ and let $r$ be its Jacobson radical. The algebra $\Lambda/r$ is semisimple. If it is separable, the Wedderburn Malcev theorem asserts that there is a subalgebra $E$ of $\Lambda$ such that $\Lambda = E \oplus r$. Note that if we assume that $k$ is algebraically closed, $\Lambda/r$ is always separable.

Consider the $E$-bimodule $r/r^2$ and the corresponding tensor algebra (see section 4 for a definition)

$$T_E(r/r^2) = E \oplus r/r^2 \oplus (r/r^2 \otimes_E r/r^2) \oplus \ldots$$

Using the evident universal property of tensor algebras, we obtain an algebra morphism

$$\varphi_\sigma : T_E(r/r^2) \to \Lambda$$

for each $E$-bimodule section $\sigma$ of the canonical surjection $r \to r/r^2$. There is always such a section since $E \otimes E^{op}$ is a semisimple algebra (recall that $E$ is separable). The algebra morphism $\varphi_\sigma$ is determined by the inclusion $E \subset \Lambda$ at zero degree and by $\sigma$ at degree one. One can show that $\varphi_\sigma$ is surjective and that $\operatorname{Ker} \varphi_\sigma$ is contained in the two-sided ideal of $T_E(r/r^2)$ generated by the degree two elements.

Different sections gives rise in general to different morphisms and different kernels. As usual, this sort of presentation of an algebra which takes into account a given maximal semisimple subalgebra is not unique.

The link to the categorical side is clear, note that if $Q$ is a finite quiver the Gabriel algebra $C_Q$ is the tensor algebra $T_{kQ_0}(kQ_1)$, where $kQ_0$ is the semisimple algebra provided by the vertices and $kQ_1$ the $kQ_0$-bimodule given by the linear combinations of arrows.

## 3 Pre-generated ideals

In this section we present part of the results of [4] which computes the dimension of $H^1$ for certain quotients of path algebras of a quiver. We begin with a definition which clarifies the situation considered in [4].

**Definition 3.1** *Let $Q$ be a finite quiver, $k$ be a field, $kQ$ the path algebra and $F$ the two-sided ideal generated by the arrows of $Q$.*

*A two-sided ideal $I$ of $kQ$ is <u>pre-generated</u> if:*

1. *$I$ is admissible (i.e. $I \subset F^2$ and there exist some $n$ such that $F^n \subset I$).*

2. Let $x$ and $y$ be two vertices of $Q$. The following alternative holds: either $yIx = y(kQ)x$ or $yIx = y(FI + IF)x$.

We provide a collection of examples of pre-generated ideals.

a) The zero ideal is pre-generated whenever $Q$ has no oriented cycles.

b) Consider an admissible two-sided ideal which is either "full or zero" between each couple of vertices $x$ and $y$. More precisely, either $yIx = ykQx$ or $yIx = 0$. Then $I$ is pre-generated.

c) A quiver is said to be narrow if for any couple of vertices $x$ and $y$ there is at most one path from $x$ to $y$. Note that each tree is narrow. Every admissible ideal of a narrow quiver is "full or zero" and hence is always pre-generated.

d) Let $Q$ be an arbitrary finite quiver, $m \geq 2$ an integer and $F^m$ the $m$-truncating ideal of paths of length at least $m$. It is easy to see that $F^m$ is a pre-generated ideal of $kQ$ if and only if all the parallel paths to each $m$-length path are of length greater or equal than $m$. For $m = 2$ the entire cohomology algebra of $kQ/F^2$ is known, see [7].

e) Let $Q$ be a $n$-cycle, namely $Q$ has $n$ vertices indexed by the cyclic group of order $n$ and linked by $n$ arrows, one from each vertex to the following one, given by the increment by a chosen generator of the group. The $m$-truncating ideal $F^m$ is pre-generated if and only if $m < n$. Note that $F^m$ is not a "full or zero" ideal.

**Proposition 3.2** *[4, p. 647]. Let $Q$ be a finite quiver, $I$ be a pre-generated ideal and $\Lambda = kQ/I$. Then*

$$\dim_k H^1(\Lambda, \Lambda) = \dim_k \mathbf{Z}(\Lambda) - \sum_{x \in Q_0} \dim_k(x\Lambda x) + \sum_{x,y \in Q_0} |yQ_1 x| \dim_k(y\Lambda x)$$

*where $yQ_1 x$ is the set of arrows from $x$ to $y$.*

**Proof**: Let $r$ be the Jacobson radical of $\Lambda$. The following is a direct consequence of the Theorem [4, p. 645]. Since $I$ is pre-generated, $H^2(\Lambda, \Lambda) = H^2(\Lambda, r^i) = 0$ for any positive integer $i$. Proposition [4, p. 647] provides the formula.

**Example 3.3** *Let $Q$ be a connected narrow quiver and $I$ an admissible ideal. Then $\dim_k H^1(\Lambda, \Lambda) = 1 - |Q_0| + |Q_1|$.*

**Remark 3.4** *If $I = 0$ and $Q$ is connected without oriented cycles, the above formula is*
$$\dim_k H^1(\Lambda, \Lambda) = 1 - |Q_0| + |Q//Q_1|$$

where we still denote by $Q$ the set of all the paths of the quiver $Q$, and $Q//Q_1$ the subset of the product $Q \times Q_1$ consisting on parallel paths, namely

$$Q//Q_1 = \{(\gamma, a) \mid s(\gamma) = s(a) \text{ and } t(\gamma) = t(a)\}.$$

Indeed, $Q$ has no cycles hence $\dim_k x(kQ)x = 1$. Moreover, $\dim_k y(kQ)x$ is the number of paths from $x$ to $y$. Note that Hochschild cohomology is additive on product of algebras (see [10]), hence if $Q$ is not connected just add the formulas for each connected component.

**Example 3.5** *Let $Q$ be the quiver with two vertices $x$ and $y$ and $n$ arrows from $x$ to $y$.*
$$\dim_k H^1(kQ, kQ) = 1 - 2 + n^2 = n^2 - 1.$$

**Example 3.6** *Let $Q$ be a $n$-cycle. Consider $I = F^m$ the truncating pre-generated ideal with $m < n$ (see examples above).*
*Then*
$$\dim_k H^1(\Lambda, \Lambda) = 1 - m + m = 1.$$

D. Happel considers in [17, p. 113] a family of algebras which are "monomial, Schurian and semi-commutative". Next we show that these algebras correspond to examples of pre-generated ideals. Actually the general monomial case will be considered in the next section.

Let $Q$ be a connected quiver without oriented cycles and let $I$ be a monomial two-sided ideal, which means that $I$ is generated by paths of length at least two. The Schurian condition on $I$ is equivalent to the fact that for each couple of vertices $(x, y)$ we have either $y(kQ)x = yIx$ or $yIx$ is one co-dimensional in $y(kQ)x$. The semi-commutative hypothesis means that if $yIx$ is one co-dimensional, then there is only one path from $x$ to $y$ in $Q$, hence $yIx = 0$. In other words this case is an example of a "full or zero" ideal, which is automatically a pre-generated ideal. The previous Proposition applies. Note that $\dim \mathbf{Z}(kQ/I) = 1$ and that each arrow has no parallel paths other than itself. Hence
$$\dim_k H^1(\Lambda, \Lambda) = 1 - |Q_0| + |Q_1|.$$

## 4 Tensor algebras

Let $E$ be a $k$-algebra and $M$ an $E$-bimodule. The tensor algebra

$$T_E(M) = E \oplus M \oplus (M \otimes_E M) \oplus \ldots$$

has a multiplication provided by the tensor product of tensors. Any algebra map $\varphi : T_E(M) \to A$ is completely determined by an algebra map $\varphi_E : E \to A$ (which provides $A$ with an $E$-bimodule structure) and an $E$-bimodule map $\varphi_E : M \to A$. This universal property characterizes $T_E(M)$.

**Example 4.1** Let $Q$ be a finite quiver, $kQ_0$ be the semisimple commutative algebra given by its vertices and $kQ_1$ the $kQ_0$-bimodule of $k$-linear combinations of arrows. Then $T_{kQ_0}(kQ_1)$ is the path algebra of $Q$ since a basis of an $m$-tensor power of $kQ_1$ over $kQ_0$ is given by the paths of length $m$ of $Q$.

**Lemma 4.2** Let $E$ be any $k$-algebra and $M$ an $E$-bimodule. There is an exact sequence
$$0 \to T \otimes_E M \otimes_E T \to T \otimes_E T \to T \to 0.$$

For a proof see [5, p. 98], [8] or [18].

**Remark 4.3** If $E$ is a separable $k$-algebra, $T \otimes_E M \otimes_E T$ is a projective $T$-bimodule. Indeed, any $E$-bimodule is projective, hence $M$ is projective. Tensoring over $E$ by a larger algebra $T$ transfers this fact to $T \otimes_E M \otimes_E T$, it is a projective $T$-bimodule. Moreover, $T \otimes_E T = T \otimes_E E \otimes_E T$ is also a projective $T$-bimodule by the same argument.

**Theorem 4.4** Let $E$ be a separable $k$-algebra, $M$ be a finitely generated $E$-bimodule, $T$ the tensor algebra and $X$ a $T$-bimodule which is finitely generated as an $E$-bimodule. Then $H^1(T, X)$ is a finite dimensional vector space and
$$\dim_k H^1(T, X) = \dim_k \mathsf{Hom}_{E-E}(M, X) - \dim_k(X^E / X^T).$$
We have $H^i(T, X) = 0$ for $i \geq 2$.

**Proof**: The above short exact sequence is a projective resolution of $T$ as a $T$-bimodule. Hence the complex of cochains
$$0 \to \mathsf{Hom}_{T-T}(T \otimes_E T, X) \to \mathsf{Hom}_{T-T}(T \otimes_E M \otimes_E T, X) \to 0$$
gives the Hochschild cohomology $H^i(T, X)$.

There are canonical isomorphisms which transforms this complex into the following:
$$0 \to X^E \to \mathsf{Hom}_{E-E}(M, X) \to 0$$
Moreover we know that $H^0(T, X) = \mathsf{Hom}_{T-T}(T, X) = X^T$.

**Remark 4.5** Since $E$ is separable, $E \otimes E^{op}$ is a semisimple algebra hence the isomorphism class of an $E$-bimodule is determined by integers (which give the multiplicity of each simple $E$-bimodule). Let $E_0$ be a complete set of central orthogonal idempotents of $E$. The isomorphism classes of simple $E$-bimodules are parameterized by $E_0 \times E_0$.

Let $({}_y n_x(M))_{x,y \in E_0}$ and $({}_y n_x(X))_{x,y \in E_0}$ be the integers which characterize $M$ and $X$. Assume the endomorphism algebra of each simple bimodule is just $k$. From the Schur Lemma we infer:
$$\dim_k \mathsf{Hom}_{E-E}(M, X) = \sum_{x,y \in Q_0} {}_y n_x(M) \, {}_y n_x(X)$$

**Corollary 4.6** Let $Q$ be a finite quiver and $X$ a $kQ$-bimodule which is finitely generated as a $kQ_0$-bimodule (i.e. $X$ is finite dimensional over $k$). Then

$$\dim_k H^1(kQ, X) = \dim_k(X^{kQ}) - \dim_k(X^{kQ_0}) + \sum_{x,y \in E_0} |yQ_1x| \dim_k(yXx)$$

where $yQ_1x$ is the set of arrows from $x$ to $y$.

**Proof**: The algebra $E = kQ_0$ is semisimple and commutative, as well as $kQ_0 \otimes (kQ_0)^{op}$. Moreover the simple bimodules are one dimensional and the multiplicity of a simple bimodule coincides with the dimension of the corresponding isotypic component. Consider $E_0 = Q_0$. We have the following:

$$_y n_x(kQ_1) = |yQ_1x|$$

$$_y n_x(X) = \dim_k(yXx).$$

**Remark 4.7** $kQ$ is finitely generated as a $kQ_0$-bimodule if and only if $Q$ has no oriented cycles.

**Corollary 4.8** Let $Q$ be a connected finite quiver without oriented cycles. Then

$$\dim_k H^1(kQ, kQ) = \dim \mathbf{Z}(kQ) - \sum_{x \in Q_0} \dim |xQx| + \sum_{x,y \in Q_0} |yQ_1x|\, |yQx|$$

where $yQx$ is the set of paths from $x$ to $y$ in the quiver.

**Corollary 4.9** Let $Q$ be any finite quiver and let $kQ/I$ be a finite dimensional quotient by a two sided ideal $I$ of $kQ$.
Then $\dim_k H^1(kQ, kQ/I) =$

$$\dim_k \mathbf{Z}(kQ/I) - \sum_{x \in Q_0} \dim_k x(kQ/I)x + \sum_{x,y \in Q_0} |yQ_1x| \dim_k y(kQ/I)x.$$

If $Q$ is connected without oriented cycles and $I$ is admissible, the formula above becomes

$$\dim_k H^1(kQ, kQ/I) = 1 - |Q_0| + \sum_{x,y \in Q_0} |yQ_1x| \dim_k y(kQ/I)x.$$

## 5 Monomial algebras

The results of this section are consequences of [5]. We present here the case where the quiver has no oriented cycles. As quoted in the Introduction, the general case is studied in a paper in with Manuel Saorin [9].

Let $Q$ be a finite quiver and $Z$ be a minimal set of paths of length at least two, which means that if $\delta$ is a path of $Z$ there is no strict sub-path of $\delta$ which belong to $Z$. The two-sided ideal $\langle Z \rangle$ generated by $Z$ has a basis $\bar{Z}$ which consists of paths containing some path of $Z$. Its complement $B$ are paths which to not contain any path of $Z$. This set $B$ is a basis of the monomial algebra $kQ/\langle Z \rangle$.

Following [5] we will define $(Q_1//B)_e$ a subset of effective couples of $Q_1//B$. First we consider glued couples: a couple $(a, \varepsilon)$ of $Q_1//B$ is glued if either $a$ is the first or the last arrow of $\varepsilon$, or if $a$ is a loop and $\varepsilon$ is the vertex of the loop. The trivial glued couples are the diagonal ones $(a, a)$ for each arrow $a$. If $Q$ has no oriented cycles, there are no glued couples besides the trivial ones.

A couple $(a, \varepsilon)$ is effective if

1. it is not glued,

2. there is a path $\gamma$ in $Z$ which contains $a$, such that the replacement of $a$ by $\varepsilon$ inside $\gamma$ provides a path which belong to $B$.

We denote $(Q_1//B)_e$ the set of effective elements and $(Q_1//B)_{ne}$ its complement. From now on, the quivers are without oriented cycles.

**Theorem 5.1** *Let $Q$ be a finite connected quiver without oriented cycles. Let $Z$ be a minimal set of paths of length at least $2$ and let $B$ be the set of paths which do not contain any path of $Z$. Then*

$$\dim_k H^1(kQ/\langle Z \rangle, kQ/\langle Z \rangle) = 1 - |Q_0| + |(Q_1//B)_{ne}|$$

**Proof**: There is an exact sequence (see [5, p. 98])

$$0 \to H^1(\Lambda, X) \to H^1(kQ, X) \to \mathsf{Hom}_{kQ-kQ}(I/I^2, X) \to H^2(\Lambda, X) \to 0$$

for $\Lambda = kQ/I$ and $X$ a $\Lambda$-bimodule.

Under the hypothesis of the Theorem we obtain:

$$\begin{aligned} 0 &\to H^1(kQ/\langle Z \rangle, kQ/\langle Z \rangle) \to H^1(kQ, kQ/\langle Z \rangle) \\ &\to \mathsf{Hom}_{kQ-kQ}(\langle Z \rangle/\langle Z \rangle^2, kQ/\langle Z \rangle) \to H^2(kQ/\langle Z \rangle), kQ/\langle Z \rangle) \to 0. \end{aligned}$$

The dimension of each term of this exact sequence is known, except the first one. Using 4.9 For the second term we have:

$$\dim_k H^1(kQ, kQ/\langle Z \rangle) = 1 - |Q_0| + \sum_{x,y \in Q_0} |yQ_1x|\, |yBx|$$

This formula becomes:

$$\dim_k H^1(kQ, kQ/\langle Z \rangle) = 1 - |Q_0| + |Q_1//B|.$$

Concerning the third term, its dimension is the number of medals, where a medal is a certain connected component of a specific quiver associated to $Q$ and $Z$, the parallel quiver (see [5]). The dimension of the last term is the number of medals minus $|(Q_1//B)_e|$ (see [5, 3.12]), the result follows. Notice that the precise definition of the parallel quiver and of a medal is unnecessary for our purposes in this paper.

Note than an element $(a, \varepsilon)$ of $Q_1//B$ is non effective if either it is glued or in case it is non-glued any $\gamma \in Z$ containing the arrow $a$ provides a path of $\bar{Z}$ when we replace $a$ by $\varepsilon$ inside $\gamma$. The diagonal $\Delta(Q_1//Q_1) = \{(a,a) \mid a \in Q_1\}$ of trivially glued couples in $Q_1//B$ is a subset of $(Q_1//B)_{ne}$. Hence, from the exact formula of the Theorem we derive the following bound:

**Corollary 5.2** *Let $Q$ be a connected quiver without oriented cycles and $Z$ be a minimal set of paths of length at least two. Then*

$$\dim_k H^1(kQ/\langle Z \rangle, kQ/\langle Z \rangle) \leq 1 - |Q_0| + |Q_1|.$$

## 6 Truncated algebras

Truncated algebras are special monomial algebras, we will apply to this algebras the results of the preceding section.

The $m$-truncating ideal $F^m$ of a quiver algebra $kQ$ is the two-sided ideal generated by $Q_m$, the set of paths of length $m$. We already stated the precise condition saying that $F^m$ is pre-generated, and we mentioned that the dimension of $H^1$ for an $m$-truncated path algebra of a $n$-cycle is 1, whenever $m < n$ (the later is the precise condition for $F^m$ to be pre-generated).

The next result gives the dimension of $H^1$ for a truncated path algebra $kQ/F^m$ ($m \geq 2$) of a quiver without oriented cycles.

**Theorem 6.1** *Let $Q$ be a finite connected quiver without oriented cycles, $\Lambda = kQ/F^m$ the $m$-truncated path algebra ($m \geq 2$). Then*

$$\dim_k H^1(\Lambda, \Lambda) = 1 - |Q_0| + |Q_1//B|.$$

**Remark 6.2** Note that $B$ is the set of paths of length strictly less than $m$.

**Proof**: Recall that $\Lambda$ is a monomial algebra with minimal set of paths $Z = Q_m$. From Theorem 5.1 we have:

$$\dim_k H^1(\Lambda, \Lambda) = 1 - |Q_0| + |(Q_1//B)_{ne}|.$$

We show that $(Q_1//B)_{ne} = Q_1//B$ in the truncated case. Recall that $(a, \varepsilon)$ is non effective if it is either glued or, in case it is non glued, any replacement of $a$ by $\varepsilon$ in a path of $Z$ provides a path of $\bar{Z}$. Note that $\varepsilon$ is not a vertex. Since $Z = Q_m$, we infer that any element of $Q_1//B$ is non effective.

The case of truncated algebras whose quiver has oriented cycles is considered in a paper with M. Saorin [9].

# 7 Finite simplicial complexes and posets

Let $S$ be a finite partially ordered set, that is a finite set equipped with an order relation $\leq$ which is reflexive, anti-symmetric and transitive. The usual incidence algebra is the vector space of functions on the subset of $S \times S$ of comparable couples $\{(y,x) \in S \times S \mid y \leq x\}$.

Let $f$ be such a function, for simplicity we denote $_yf_x$ the scalar $f(y,x)$ for $y \leq x$. The product is given by

$$_y(fg)_x = \sum_{y \leq z \leq x} {_yg_z}\, {_zf_x}.$$

Considering elementary functions (Dirac masses on pairs $(y,x)$ such that $y \leq x$) it is clear that the incidence algebra is isomorphic to an admissible quotient of a quiver algebra: let $Q$ be the quiver which has as vertices the set $S$ and which has an arrow from $x$ to $y$ in case $x \geq y$ with no proper $z$ between $x$ and $y$. Consider the ideal of $kQ$ generated by all the differences of parallel paths. This ideal is admissible and $kQ/I$ is canonically isomorphic to the incidence algebra.

Recall that a simplicial complex $\Sigma$ is a set equipped with an "hereditary" non empty collection of subsets, the simplices: if $\sigma$ is a simplex of the collection, every subset of $\sigma$ is a simplex. In other words, the simplicial complex is completely determined by a non empty collection of maximal subsets of $\Sigma$. Of course, a simplicial complex provides a poset: the elements are the simplices and the order relation is the inclusion (called incidence relation).

Conversely, a poset $S$ gives a simplicial complex on the set $S$, the simplices are the chains of ordered elements of $S$.

Actually the simplicial complex associated to the poset of a simplicial complex $\Sigma$ is the "barycentric subdivision" of $\Sigma$. Gerstenhaber and Schack proved that the Hochschild cohomology of the incidence algebra and the usual simplicial cohomology are canonically isomorphic ([15], [16]). See also [6].

**Example 7.1** *Consider the poset $S = \{a,b,c,d\}$ with order relation given by $a > c$, $a > d$, $b < c$, $b < d$. Its quiver $Q$ has four vertices and four arrows, the ideal for the incidence algebra is 0. We already know that*

$$\dim_k H^1(kQ, kQ) = 1 - 4 + 4 = 1.$$

*Note that the simplicial complex associated to $S$ has $S^1$ as geometric realization, we also know that $\dim_k H^1(S^1, k) = 1$.*

Département de Mathématiques, Université de Montpellier 2,
F-34095 Montpellier cedex 5.
cibils@math.univ-montp2.fr
January 1999